\documentclass[preprint,12pt]{elsarticle} 
\usepackage{bm}
\usepackage{amsmath,amsfonts,amssymb,amsthm,mathtools} 
\usepackage{hyperref}
\usepackage{listings}
\usepackage{epstopdf}
\usepackage{algorithm}
\usepackage[]{algpseudocode}
\usepackage{bbm}
\makeatletter
\def\BState{\State\hskip-\ALG@thistlm}
\makeatother

\usepackage{amssymb}
\usepackage{todonotes}
\usepackage{amsthm}
\newtheorem{proposition}{Proposition}

\newtheorem{lemma}{Lemma}

\newcounter{bla}

\journal{Computer Physics Communications}

\begin{document}

\begin{frontmatter}



\title{Numerical solution of the Boltzmann equation with S-model collision integral using tensor decompositions}


\author[a]{A.V. Chikitkin \corref{author}}
\author[a]{E.K. Kornev}
\author[a,b]{V.A. Titarev}

\cortext[author] {Corresponding author.\\\textit{E-mail address:} chikitkin.av@mipt.ru}

\address[a]{Moscow Institute of Physics and Technology, 9 Institutskiy per., Dolgoprudny, Moscow Region, 141701, Russian Federation}

\address[b]{Federal Research Center “Computer Science and Control”, Russian Academy of Sciences, Vavilov str. 40, Moscow, 119333, Russian Federation}

\begin{abstract}

The paper presents a new solver for the numerical solution of the  Boltzmann kinetic equation with the Shakhov model collision integral (S-model) for arbitrary spatial domains. The numerical method utilizes the Tucker decomposition, which reduces the required computer memory for up to 100 times, even on a moderate velocity grid. 
This improvement is achieved by representing the distribution function values on a structured velocity grid as a 3D tensor in the Tucker format. The resulting numerical method makes it possible to solve complex 3D problems on modern desktop computers. Our implementation may serve as a prototype code for researchers concerned with the numerical solution of kinetic  equations in 3D domains using a discrete velocity method.





\end{abstract}

\begin{keyword}
rarefied gas dynamics \sep Shakhov model \sep discrete velocity method \sep unstructured mesh \sep Tucker decomposition \sep tensor decompositions \sep finite-volume methods

\end{keyword}

\end{frontmatter}



{\bf PROGRAM SUMMARY}

\begin{small}
\noindent
{\em Program Title:} Boltzmann-T                                      \\
{\em Licensing provisions:} MIT                                  \\
{\em Programming language:} Python 3     \\
{\em External libraries:} Solver is based on the customized version of the tucker3d  library  [1] \\
{\em Nature of problem:} Numerical solution of the Boltzmann kinetic equation with the S-model collision integral in an arbitrary 3D spatial  domain \\
{\em Solution method:} Discrete velocity method utilizing tensor decomposition for memory reduction \\
{\em Restrictions:} At present, 1st order advection scheme is used, solver supports unstructured hexagonal meshes written in StarCD ASCII format \\
{\em Additional comments:} Source code is available at \url{https://github.com/chikitkin/Boltzmann-Tucker}

\end{small}

\section{Introduction}
\label{sec:intro}
The Boltzmann kinetic equation (BKE) is the main mathematical model of the theory of rarefied gases. Due to the high dimensionality of the phase space and the complexity of the collision integral, the numerical solution of the BKE is much more complicated and computationally expensive than the numerical solution of macroscopic equations, such as the Navier-Stokes equations of the compressible gas~\cite{Petrov:2018a}.

There are several simplified collision operators which preserve a number of  important properties of the exact operator. The simplest operator results in the BGK model \cite{betahatnagar1954model}. A more accurate approximation is given by the Shakhov model (S-model) \cite{Shakhov1968} and its extention to the diatomic gases by Rykov \cite{Rykov1975}. Comparisons with calculations using the exact Boltzmann equation, the direct simulation Monte Carlo method, and with the experimental data have confirmed high accuracy of the S-model, see  \cite{Sharipov:1998a,Titarev:2012d,Frolova:2018b,titarev2018application,Titarev2020} and the references therein.

In model kinetic equations, the calculation of the collision integral requires only a certain number of macroparameters or moments of the distribution function, which are computed via 3-dimensional integrals over the velocity space. Despite this simplification, numerical solution of the simplified equation is still quite a  demanding computational task, especially for three-dimensional applications. One of the approaches to reduce the computational cost and memory requirements of numerical methods for model kinetic equations is the use of an adaptive unstructured mesh in the velocity space \cite{Arslanbekov:2013a,Baranger:2014,guo:sparse,titarev2016openmp,titarev2018application}. It should be noted that it significantly complicates the numerical algorithm and often requires some a-priori information about the problem to be solved. The simplest algorithm can be constructed for structured Cartesian grids in the velocity space. In this case,  the distribution function values at all nodes of the mesh form a multidimensional array, which will be hereafter called ``tensor''. Therefore, the natural way to speed up the method and reduce the required memory is to employ low-rank tensor approximations, which are well-known in linear algebra. This is justified by theoretical estimates proving that for tensors, generated by the values of smooth functions, such approximations always exist \cite{Tyrtyshnikov2004423,Tyrtyshnikov2004Sbornik}.

There are many studies on this subject. In \cite{Khoromskij20071291}, a special tensor format is proposed for the approximation of tensors which arise from calculating the exact collision integral on a tensor grid. In \cite{Dolgov2014268}, tensor decompositions were successfully applied to the numerical method for the Vlasov equation with the BGK collision integral. The memory consumption was reduced 17 times as compared with the standard numerical method on the same meshes.
Another version of the numerical method for the Vlasov equation is described in \cite {Kormann2015B613}. It is noted that the use of tensor decompositions reduces storage by more than 100 times.
A recent paper \cite{Boelens2018519} describes the general framework for the application of tensor decompositions to the numerical methods for partial differential equations of a certain type. The results of test calculations of simple problems for the BGK model collision integral are presented.

In the cited papers, tensor decompositions are applied to tensors formed by values of the distribution function on a structured tensor grid in both physical and velocity space. Such tensors have dimension 4 or 6 depending on the dimensionality of the problem. For such dimensions, low-rank approximations are especially effective. However, this approach is applicable only to problems with a simple shape of the computational domain and simple boundary conditions so that one can use a structured grid in the physical space. On the contrary, in many applications, the computational domain has a complex shape. For such problems, an unstructured mesh in the physical space is often used. In this regard, it is more convenient to approximate a tensor formed by the distribution function values only on the velocity grid at each point of the physical space.

In this paper, we propose an analog of the discrete velocity method in which tensors formed by the distribution function values on the velocity grid are approximated using the Tucker format \cite{tucker1963}. The examples of test calculations are presented, which show that the proposed approach allows reducing the computer memory consumption by 100 times while maintaining satisfactory accuracy; the CPU time increases only mildly.

\section{Mathematical model}
\label{sec:model}
The Boltzmann equation of a monatomic gas with a model collision integral has the following form:
\begin{align}
    & \frac{\partial f}{\partial t} + \sum_{\alpha = 1}^{3} \xi_{\alpha} \frac{\partial f}{\partial x_{\alpha}} = J(f),
    \label{eq:model-eq}
\end{align}
where $f(t, \bm{x}, \bm{\xi})$ is the value of the distribution function, $\bm x$ is the space coordinates, $\bm{\xi}$ is the velocity vector, $J$ is the collision operator. In the Shakhov model \cite{Shakhov1968}, the collision operator is defined by the vector of macroparameters $\bm a = \bm a (f) = [n, \bm u, T, \bm q]$, which are expressed through the moments of the distribution function: 
\begin{align}
\begin{split}
    & n = \int fd\bm{\xi} , \quad \bm{u} = \frac{1}{n} \int \bm{\xi} f d\bm{\xi}, \\
    & T = \frac{1}{3 n R_g} \left(\int \xi^2 f \, d\bm \xi - n u^2 \right), \quad \rho = mn, \;  p = \rho R_{g}T \\
    & \bm{v} = \bm{\xi} - \bm{u}, \quad \bm q = \frac{1}{2} m \int \bm v v^2 f \, d \bm \xi .
\end{split}
\end{align}
Here, $n$ is the numerical density, $\bm u$ is the macro velocity, $R_g$ is the gas constant, $T$ is the temperature, $\rho$ is the mass density, $p$ is the pressure, $m$ is the mass of one molecule, $\bm q$ is the heat flux vector. All the macroparameters are functions of $t$ and $\bm x$. 
The expression for the function $J(f)$ is the following:
\begin{align}
\begin{split}
    & J(f)(t, \bm x, \bm \xi) = \nu(t, \bm x) (f^S(t, \bm x, \bm \xi; \bm a) - f(t, \bm x, \bm \xi)),\; \nu = \frac{p}{\mu(T)} \\
    & f^S(t, \bm x, \bm \xi ;\bm a) = f^{M} \bigg[1 +  \frac{4}{5}(1 - Pr) \frac{n}{m (2 R_g T)^2}\left(\sum_{\alpha = 1}^3 q_\alpha v_\alpha\right) \left(\frac{v^2}{2 R_g T} - \frac{5}{2} \right) \bigg] \\
    & f^{M} = f^{M}(t, \bm x, \bm \xi; \bm a) = \frac{n}{(2\pi R_{g}T)^{3/2}}\exp \left(-\frac{v^2}{2 R_g T} \right) .
\end{split}
\label{eq:shakhov_integral}
\end{align}
Here $\mu = \mu(T)$ is the dynamic viscosity, $Pr = 2/3$ is the Prandtl number, $f^M$ is the Maxwell (equilibrium) distribution function for the macroparameters~$\bm a$.

At the boundaries of the computational domain in the physical space, it is necessary to specify distribution function values for molecules whose velocity vector is directed inside the domain. On the surface of the body, the boundary condition of diffuse reflection with full thermal accommodation to the surface temperature $T_w$ is used. The distribution function of reflected molecules is written as:
\begin{align}
    f_w(\bm \xi) = \frac{n_w}{(2 \pi R_g T_w)^{3/2}} \exp\left(-\frac{\xi^2}{2 R_g T_w}\right) .
    \label{eq:wall-bc-f}
\end{align}
The density $n_w$ of the reflected molecules is found from the impermeability condition:
\begin{align}
    \int\limits_{\xi_n > 0} \xi_n f \, d \bm{\xi} + \int\limits_{\xi_n < 0} \xi_n f_w \, d \bm{\xi} = 0,
\end{align}
where $\xi_n = \bm \xi \cdot \bm e$ is the projection of the velocity onto the unit normal to the surface $\bm e$, directed outside the computational domain, $f$ is the distribution function of molecules coming to the wall.

For a plane of symmetry, the following boundary condition is set:
\begin{align}
    f(t, \bm{x}, \bm{\xi}) = f(t, \bm{x}, \tilde{\bm{\xi}}), \; \tilde{\bm{\xi}} = \bm{\xi} - 2 (\bm \xi \cdot \bm e) \bm e,
\end{align}
where $\bm e$ is the outward looking unit normal vector for the plane of symmetry.

For the free stream boundary condition, the distribution function at the boundary is equal to the Maxwell distribution function for the prescribed values of macroparameters.
%
\section{Discrete velocity method}
\label{sec:discr-vel}
In this paper, we use a variant of the discrete velocity method described in \cite{titarev2010jvm}, \cite{titarev2016openmp}, \cite{Titarev201417}. For brevity, we explain the main idea using an explicit first-order method, although a scheme of arbitrary approximation order can be used. Furthermore, some implicit schemes can also be used as will be shown below. 

Without loss of generality, we consider a uniform Cartesian grid in the velocity space with the following nodes:
\[\xi_{\alpha, i_\alpha} = \xi_{min} + (i_\alpha-1)\Delta \xi, \; i_\alpha = 1, \ldots, N, \; \alpha = 1,2,3, \]
where $\Delta \xi$ is the step in the velocity grid.
In general, the one dimensional grids in each dimension can be different and nonuniform. The integrals in the velocity space for any function $g$ are approximated by the 2nd order quadrature formula:
\begin{align}
    \int g(\bm{\xi}) \, d \bm{\xi} \approx (\Delta \xi)^3 \sum_{i_1, i_2, i_3 = 1}^N g(\xi_{1, i_1}, \xi_{2, i_2}, \xi_{3, i_3}) .
\end{align}
The distribution function values at the nodes of the velocity grid form a three-dimensional tensor, which is denoted by $\hat{f}(t, \bm{x})$:
\begin{align}
    \hat{f}(t, \bm{x})(i_1, i_2, i_3) = f(t, \bm{x}, \xi_{1, i_1}, \xi_{2, i_2}, \xi_{3, i_3}), \; i_1, i_2, i_3 = 1, \ldots, N .
\end{align}
Here, the second round brackets after $\hat{f}$ contain indices of the tensor $\hat{f}(t, \bm{x})$. We will further denote all the tensors by a symbol with hat to distinguish them from scalars.
Writing the kinetic equation at each node of the velocity grid, we obtain a system of $N^3$ linear constant-coefficient equations with a source term. This system can be written in the tensor form:
\begin{align}
    \hat{f}_t + (\hat{\xi}_1 \circ \hat{f})_{x_1} + (\hat{\xi}_2 \circ \hat{f})_{x_2} + (\hat{\xi}_3 \circ \hat{f})_{x_3} = \nu (\hat{f}^S - \hat{f}) .
    \label{eq:tensor-pde}
\end{align}
Here, $\hat{\xi}_{\alpha}$ is the tensor formed by the values of the corresponding velocity component at each node of the velocity grid:~$\xi_\alpha(i_1, i_2, i_3) = \xi_{\alpha, i_\alpha}$; ``$\circ$'' denotes the component-wise (Hadamard) product of tensors.

A standard finite-volume method of the Godunov type is used to discretize the left-hand side of the system \eqref{eq:tensor-pde}. The computational domain in the physical space is divided into $N_C$ finite volumes (polyhedrons)  $V_i$, $ i = 1, \ldots, N_C $. System \eqref{eq:tensor-pde} is integrated over $V_i$, the volume integral is replaced by the sum of surface integrals over the cell faces from the fluxes projected onto the normal to the face. We assume that all the faces in the mesh are numbered with the index $j = 1, \ldots, N_F$. Thus we obtain  a semi-discrete scheme of the following form:
\begin{align}
    \begin{split}
    & \frac{d \hat{f}_i}{d t} = -\frac{1}{\lvert V_i \rvert} \sum_{l = 1}^{N_F(i)} \hat{\Phi}_{j(i,l)} \, \mbox{sign}(i, l) + \hat{J}(\hat{f}_i), \; i = 1, \ldots, N_C \\
    & \hat{\Phi}_j = \int_{A_j} \hat{\xi}_{n,j} \circ \hat{f}(t, \bm x) \, dS, \; \hat{\xi}_{n,j} = e_{j, 1} \, \hat{\xi}_1 +  e_{j, 2} \, \hat{\xi}_2 +  e_{j,3} \, \hat{\xi}_3 .
    \end{split}
    \label{eq:ode_system}
\end{align}
Here, $\hat{f}_i = \hat{f}_i(t)$ is the tensor consisting of the integral averages over cell $V_i$ of the distribution function values, $|V_i|$ is the volume of the cell, $N_F(i)$ is the number of faces in the cell with index $i$, $j(i, l)$ is the global index of the $l$-th face of the cell with index $i$, $\mbox{sign}(i,l)$ equals $+1$ if the normal vector for the face $j(i,l)$ is outer with respect to the cell $i$ and $-1$ otherwise, $\bm{e}_{j}$ is the unit normal vector of the face with global index $j$, $\hat{\xi}_{n,j}$ is the projection of the velocity onto the normal vector, $A_j$ is the face with index $j$, $\hat{J}(\hat{f}_i)$  is the tensor with elements computed by formulas \eqref{eq:shakhov_integral} in which all the integrals are replaced with the quadrature formula. Face flux $\hat{\Phi}_{j}$ is approximated via the solution of the one-dimensional Riemann problem:
\begin{align}
    \hat{\Phi}_j \approx |A_j| \hat{F}(\hat{f}_{L, j}, \hat{f}_{R, j}) ,
\end{align}
where $|A_j|$ is the area of the face, $\hat{F}$ is the exact or approximate solution (flux) to the Riemann problem, $\hat{f}_{L, j}$ and $\hat{f}_{R, j}$ are the reconstructed values to the left and to the right of the face with respect to the normal vector. For the first order reconstruction at the inner faces, these values are simply values of the integral average in the left and right cell, respectively. If the face lies on the boundary then one of these values is set based on the corresponding boundary condition. The final form of the numerical method depends on the specific flux approximation and the time-marching scheme.

In the case of the exact solution to the Riemann problem (the CIR scheme \cite{toro2013riemann}), the expression for the face flux is the following:
\begin{align}
    \hat{F}_j = \hat{F}(\hat{f}_{L, j}, \hat{f}_{R, j}) = \frac{1}{2} \hat{\xi}_{n,j} \circ \left(\hat{f}_{L, j} + \hat{f}_{R,j} \right) - \frac{1}{2} \lvert \hat{\xi}_{n,j} \rvert \circ \left( \hat{f}_{R, j}  - \hat{f}_{L, j}\right) .
    \label{eq:exact-flux}
\end{align}
It should be noted that in \eqref{eq:exact-flux} instead of $\lvert \hat{\xi}_{n,j} \rvert$ some estimate can be used. This may be interpreted as using a Riemann solver of the Rusanov type \cite{toro2013riemann}. Another important fact is that the term $\lvert \hat{\xi}_{n,j} \rvert$ (or its approximation) is responsible for the level of the numerical dissipation of the scheme.

Using the explicit Euler method to solve the ODE system \eqref{eq:ode_system}, we obtain the fully discrete method:
\begin{align}
\begin{split}
    & \frac{\hat{f}^{n+1}_i - \hat{f}^n_i}{\Delta t} = \hat{R}_i^n =  -\frac{1}{\lvert V_i \rvert} \sum_{l = 1}^{N_F(i)} \mbox{sign}(i,l) |A_{j(i,l)}| \hat{F}_{j(i,l)}  + \hat{J}_i (\hat{f}_i^n), \\
    & i = 1, \ldots, N_C ,
\end{split}
\label{eq:expl-scheme-final}
\end{align}
where $\Delta t$ is the time step, superscript $n$ is the index of the time level. We will consider only a steady test problem in which the numerical solution converges to some steady state and, therefore, $n$ can be considered as the iteration index. The main part of the procedure for performing one time step is listed in algorithm \ref{alg1}. Time level indices are omitted since actually only one set of values are stored in memory.
\begin{algorithm}[H]
\caption{Time step}
\label{alg1}
\begin{algorithmic}[1]
\State \ldots \Comment{set boundary conditions}
\For{$j = 1, N_F$} \Comment{fluxes on faces}
    \State $\hat{F}_j = \frac{1}{2} \hat{\xi}_{n,j} \circ \left(\hat{f}_{L, j} + \hat{f}_{R,j} \right) - \frac{1}{2} \lvert \hat{\xi}_{n,j} \rvert \circ \left( \hat{f}_{R, j}  - \hat{f}_{L, j}\right)$ 
\EndFor
\For{$i = 1, N_C$} \Comment{compute right-hand side}
    \State $\hat{R}_i = $ \textbf{computeJ}($\hat{f}_i$) \Comment{compute collision integral}
    \For{$l = 1, N_F(i)$} \Comment{loop over faces of the cell $i$}
    \State $\displaystyle \hat{R}_i = \hat{R}_i - \mbox{sign}(i,l) \frac{|A_{j(i,l)}|}{|V_i|} \hat{F}_{j(i,l)}$ \Comment{add flux with sign}
    \EndFor
\EndFor
\For{$i = 1, N_C$} \Comment{update values}
    \State $\hat{f}_i = \hat{f}_i + \Delta t \hat{R}_i$
\EndFor
\end{algorithmic}
\end{algorithm}
The pseudo-code of the function for computing the model collision integral is given in algorithm \ref{alg2}. The function ``sum($\cdot$)'' calculates the sum of all elements in the tensor, the symbol $\hat {\mathbbm{1}} $ denotes the tensor consisting of ones: $\hat {\mathbbm{1}}(i_1, i_2, i_3) = 1$.
\begin{algorithm}[H]
\caption{Calculation of the collision integral}
\label{alg2}
\begin{algorithmic}[1]
\Procedure{computeJ}{$\hat{f}$}
\State $n = (\Delta \xi)^3 \mbox{ sum}(\hat{f})$ 
\State $u_\alpha = (\Delta \xi)^3 \mbox{ sum}(\hat{\xi}_\alpha \circ \hat{f}) / n, \; \alpha = 1,2,3$
\State $ \displaystyle u^2 = \sum_{\alpha = 1}^3 u_\alpha^2$
\State $\displaystyle \widehat{\xi^2} = \sum_{\alpha = 1}^3 \hat{\xi}_\alpha \circ \hat{\xi}_\alpha$ 
\State $\displaystyle T = \frac{1}{3 n R_g} \left((\Delta \xi)^3 \mbox{ sum}(\widehat{\xi^2} \circ \hat{f}) - n u^2\right)$, $\;$ $\rho = m n, \; p =  \rho R_g T$
%
%
\State $\displaystyle \hat{v}_\alpha = \hat{\xi}_\alpha - u_\alpha \hat{\mathbbm{1}}, \; \alpha = 1,2,3$
\State $\displaystyle \widehat{v^2} = \sum_{\alpha = 1}^3 \hat{v}_\alpha \circ \hat{v}_\alpha$
\State $\displaystyle q_\alpha = \frac{1}{2} (\Delta \xi)^3 \mbox{ sum}(\hat{v}_\alpha \circ \widehat{v^2} \circ \hat{f}), \; \alpha = 1,2,3$
\State $\displaystyle \hat{f}^M = \frac{n}{(2 \pi R_g T)^{3/2}} \exp \left(-\frac{\widehat{v^2}}{2 R_g T} \right)$
\State \begin{small}$\displaystyle \hat{f}^S = \hat{f}^M \circ \left( \hat{\mathbbm{1}} + \frac{4}{5} \left(1 - Pr\right) \frac{n}{m (2 R_g T)^2} \left(\sum_{\alpha = 1}^{3} q_\alpha \hat{v}_\alpha \right) \circ \left( \frac{\widehat{v^2}}{2 R_g T} - \frac{5}{2} \hat{\mathbbm{1}}\right)\right)$ \end{small}
\State $\displaystyle \hat{J} = \frac{p}{\mu(T)} \left(\hat{f}^S - \hat{f} \right)$
\State \textbf{return} $\hat{J}$
\EndProcedure
\end{algorithmic}
\end{algorithm}
The main observation which can be made from the listed algorithms is that one step of the numerical method requires only a few simple operations with three-dimensional tensors, namely:
\begin{enumerate}
    \item component-wise sum of two tensors
    \item component-wise product of two tensors
    \item sum of all elements in a tensor or, in the case of a nonuniform Cartesian grid in the velocity space, the convolution of the following form: 
    \begin{align}
        S = \sum_{i_1, i_2, i_3} \hat{f}(i_1, i_2, i_3) w_1(i_1) w_2(i_2) w_3(i_3) ,
    \end{align}
    where $w_\alpha$ are the 1D vectors consisting of weights of a quadrature rule.
\end{enumerate}
It follows from this observation that if there is some parametric representation of tensors for which all these operations can be performed then the storage of all tensor elements can be avoided.

The same conclusion is true for many implicit methods. In our code, we implemented a version of the LU-SGS (Lower-Upper Symmetric-Gauss-Seidel) method. This method is very effective since its computational cost is only about 50\% larger then the cost of the explicit method. For brevity, we do not list all the formulas. The details of the implementation in the context of kinetic solvers can be found in \cite{Titarev:2012c,titarev2016openmp,Chikitkin2018503}. 

In the next section, we briefly formulate the general idea of tensor decompositions and describe the Tucker decomposition which is used for the modification of the described discrete velocity method.

\section{Tensor decompositions}
Tensor decompositions extend the idea of separation of variables to multidimensional arrays. In the two-dimensional case, for any real matrix of rank $r$ there exists the singular value decomposition (SVD):
\begin{align}
    A = U \Sigma V^T, \; A(i_1, i_2) = \sum_{k = 1}^r \sigma_k u_k(i_1) v_k(i_2) ,
    \label{eq:svd}
\end{align}
where $U$ and $V$ are unitary matrices, $u_k, v_k$ are their columns.
The Eckart-Young theorem states that the best approximation of the rank $r' < r$ to the matrix $A$ in the 2-norm and the Frobenius norm is obtained by dropping the $r-r'$ terms in the SVD of $A$, which correspond to the smallest singular numbers. The low-rank approximation allows  to reduce the required memory to $2 n r'$, where $n$ is the size of the matrix (for the case of a square matrix). In this section, we will denote by $n$ the integer size of a tensor or a matrix, not the numerical density.

A direct generalization of the form \eqref{eq:svd} in the multidimensional case is the canonical decomposition (CANDECOMP, PARAFAC) \cite{DeLathauwer20001324}:
\begin{align}
    A(i_1, \ldots, i_d) = \sum_{k = 1}^r U_1(i_1, k) \ldots U_d(i_d, k) .
\end{align}
The minimal number $r$ required to express $A$ is called the tensor rank or the canonical rank.
The application of the canonical decomposition in numerical methods is limited because for $d \ge 3$ the problem of finding an approximation with a fixed canonical rank can be ill-posed \cite{de2008tensor}. Nevertheless, there are theoretical estimates showing that  tensors formed by values of a smooth function on a structured grid can be approximated with high accuracy by a low-rank tensor \cite{Tyrtyshnikov2004Sbornik}. 

For low dimensions, the Tucker decomposition is often used \cite{tucker1963}:
\begin{align}
    A(i_1, \ldots, i_d) = \sum_{k_1, \ldots, k_d= 1}^{r_1, \ldots, r_d} G(k_1, \ldots, k_d) U_1(i_1, k_1) \ldots U_d(i_d, k_d) .
\end{align}
The tensor $G$  is called the Tucker core (we will not use hat for it) and matrices $U_\alpha$, $\alpha = 1, \ldots, d$ are referred to as Tucker factors or mode factors, $r_\alpha$ are referred to as mode ranks or Tucker ranks. Sometimes the factors are assumed to be column-wise orthonormal. For the simplicity of notation, we will further assume in all complexity estimates that $r_\alpha = r$, $n_\alpha = n$ and will omit summation limits in some formulas.

This representation allows employing robust SVD-based procedures for basic linear algebra operations with tensors. The Tucker decomposition does not circumvent the ``curse of dimensionality'' since the number of parameters is $O(r^d + d n r)$ and grows exponentially in $d$. However, in many problems, ranks are very small and the Tucker decomposition turns out to be very effective for low dimensions. In this paper, we use the Tucker decomposition for the representation of 3D tensors.

It is worth mentioning two tensor formats applicable to arbitrary dimension $d$, which generalize the idea of the Tucker format: the hierarchical-Tucker (HT) format \cite {Grasedyck2009}, and the Tensor-Train (TT) format \cite {Oseledets2011}. Both formats are based on a dimensionality reduction tree and can utilize SVDs of auxiliary matrices for a low-rank approximation of a tensor. In these formats, the amount of storage and complexities of basic algorithms do not grow exponentially in $d$.

Below, we list all operations with tensors in the Tucker format which are used to modify the baseline discrete velocity method: 
\begin{enumerate}
    \item For a ``full'' tensor $A$ compute tensor $B$ in the Tucker format which approximates $A$ with a prescribed accuracy:
    \[\Vert A - B \Vert_F \le \epsilon \Vert A \Vert_F .\]
    Here $\Vert \cdot \Vert_F$ is the Frobenius norm. This can be done using the high-order SVD algorithm from \cite{de2000multilinear} with the complexity $O(n^4)$.
    \item Compute element-wise sum $C$ of two Tucker tensors $A$ and $B$. The core and the factors of $C$ are expressed directly via cores and factors of $A$ and $B$:
    \begin{align}
    \begin{split}
        & G^C(k_1, k_2, k_3) = 
        \begin{cases}
        G^A(k_1, k_2, k_3), & \mbox{if } 1 \leq k_\alpha \leq r^A_\alpha \\
        G^B(k_1 - r^A_1, k_2 - r^A_2, k_3-r^A_3), & \mbox{if } r^A_\alpha < k_\alpha \leq r^A_\alpha + r^B_\alpha \\
        0, & \mbox{else}
        \end{cases} \\
        & k_\alpha = 1, \ldots, r^A_\alpha + r^B_\alpha \\
        & U^C_\alpha =
        \begin{bmatrix}
        U^A_\alpha & U^B_\alpha
        \end{bmatrix}.
    \end{split}
    \label{eq:element_wise_sum}
    \end{align}
    Here, the superscripts mark cores, factors and ranks of the corresponding tensor. Element-wise sum does not require any calculations but the ranks of the sum are equal to the sum of the  ranks of the summands.
    \item Compute the element-wise (Hadamard) product $C = A \circ B $. The core and the factors of the result can be computed using these formulas:
    \begin{align}
    \begin{split}
        & G^C(k_1, k_2, k_3) = G^A(k_1^A, k_2^A, k_3^A)  G^B(k_1^B, k_2^B, k_3^B), \; \\
        & k^A_\alpha = \lceil k_\alpha / r_\alpha^B \rceil, \quad k_\alpha^B =  mod(k_\alpha, r_\alpha^B) + 1, \;  k_\alpha = 1, \ldots, r^A_\alpha r^B_\alpha\\
        %
        & U^C_\alpha(i_\alpha, :) = U^A_\alpha(i_\alpha, :) \otimes U^B_\alpha(i_\alpha, :) .
    \end{split}
    \end{align}
    Here $\lceil \cdot \rceil$ is the ceiling function, $mod(a,b)$ is the remainder of the division of $a$ by $b$, $\otimes$ is the Kronecker product, and $U^A_\alpha(i_\alpha, :)$ is the row of the matrix (we use Fortran/Matlab slicing notation).

    The element-wise multiplication requires $O(r^6 + nr^2)$ operations; the ranks of the product are equal to the product of the ranks of $A$ and $B$.
    \item 
    Rounding or recompression: given a tensor in the Tucker format with the ranks $r_0$ approximate it by a tensor in the same format with the ranks $r < r_0$.  This can be done by the algorithm 3 from \cite{oseledets2009linear} with the complexity $O(n r_0^2 + n r_0 r + r_0^4)$. The final ranks can be determined by the prescribed relative accuracy $\epsilon$.
    \item For a tensor $A$ in the Tucker format compute convolution with one dimensional vectors:
    \begin{align}
        S = \sum_{i_1, i_2, i_3} A(i_1, i_2, i_3) w_1(i_1) w_2(i_2) w_3(i_3) .
    \end{align}
    Substituting the representation of $A$ we obtain:
    \begin{align}
    \begin{split}
        & S = \sum_{k_1, k_2, k_3} G(k_1, k_2, k_3) S_1(k_1) S_2(k_2) S_3(k_3), \\
        & S_\alpha(k_\alpha) = \sum_{i_\alpha} U_\alpha(i_\alpha, k_\alpha) w_\alpha(i_\alpha), \; \alpha = 1, 2, 3 .
    \end{split}
    \end{align}
    Using this formulas, convolution is computed with complexity $O(n r^3)$.
    %
\end{enumerate}
All the listed basic procedures allow rewriting the algorithm of the discrete velocity method as a sequence of operations with tensors in the Tucker format: element-wise operations are replaced by their analogs, besides, intermediate rounding is added to prevent the growth of the Tucker  ranks.
The next section describes the details of the adaptation of the algorithm.
\section{Tensorized discrete velocity method}
\label{sec:tens-discr-vel}
In the tensorized version of the method, all low-rank arrays are constructed directly in the Tucker form. 
For example, since the Maxwell distribution function is the product of three 1D functions, we can explicitly construct the Tucker tensor with ranks 1 and with the corresponding factors:
\begin{align}
    \begin{split}
        & \hat{f}^M(i_1, i_2, i_3) = U_1(i_1) U_2(i_2) U_3(i_3), \\
        & U_\alpha(i_\alpha) = \frac{n^{1/3}}{(2 \pi R_g T)^{1/2}} \exp \left(- \frac{(\xi_{\alpha, i_\alpha} - u_{\alpha})^2}{2 R_g T} \right) .
    \end{split}
\end{align}
Since the tensor for the Shakhov function is computed from simple tensors of ranks 1, we can easily prove estimates for its ranks. First, we prove an auxiliary statement:
\begin{lemma} 
For any vectors $w_1(i_1), w_2(i_2), w_3(i_3)$ tensor $\hat{A}$ of the form
\[\hat{A}(i_1, i_2, i_3) = w_1(i_1) + w_2(i_2) + w_3(i_3)\]
can be expressed in  the Tucker format with ranks $2$.
\label{lemma:lemma-1}
\end{lemma}
\begin{proof}
It can be verified by the direct computation that the Tucker tensor with the following core and factors is equal to $\hat{A}$:
\begin{align}
\begin{split}
    & G(k_1, k_2, k_3) = 0, \mbox{ except for } G(1,2,2) = G(2, 1, 2) = G(2,2,1) = 1, \\
    & U_\alpha(:, 1) = w_\alpha, \; U_\alpha(:, 2) = 1, \alpha = 1,2,3 .
\end{split}
\end{align}
\end{proof}
Using this observation, we can prove the following proposition.
\begin{proposition}
For any velocity grid, tensor $\hat{f}^S$ can be represented in the Tucker format with ranks $5$.
\end{proposition}
\begin{proof}
Expression for $\hat{f}^S$ has the following structure:
\begin{align}
    \hat{f}^S = \hat{f}^M \circ \left(\hat{\mathbbm{1}} +  \left(\sum_{\alpha = 1}^{3} c_\alpha \hat{v}_\alpha \right) \circ \left( b \widehat{v^2} - \hat{\mathbbm{1}}\right)\right)
    \label{eq:fshakhov_structure}
\end{align}

Here $c_\alpha$ and $b$ are some scalars, and $\hat{\mathbbm{1}}$ and $\hat{f}^M$ have ranks $1$. By the lemma \ref{lemma:lemma-1}, ranks of tensors $ \displaystyle \left(\sum_{\alpha = 1}^{3} c_\alpha \hat{v}_\alpha \right)$ and $\left( b \widehat{v^2} - \hat{\mathbbm{1}}\right)$ are equal to $2$, since $\hat{v}_\alpha (i_1, i_2. i_3) = \xi_{\alpha,i_\alpha} - u_\alpha$.

Therefore, the ranks of the expression in brackets in \eqref{eq:fshakhov_structure} computed using element-wise operations are equal to $1 + 2 \cdot 2 = 5$. Since the ranks of $\hat{f}^M$ equal to $1$, the ranks of the $\hat{f}^S$ are the same.
\end{proof}
The only bottleneck in the algorithm is  the tensor $|\hat{\xi}_{n,j}|$ and the tensors $\hat{\xi}_{n,j}^+$, $\hat{\xi}_{n,j}^-$, in which the negative (positive) values in the tensor $\hat{\xi}_{n,j}$ are replaced with zero. In the general case, the normal vector $\bm{e}_j$ does not coincide with one of the coordinate axes. Then the mentioned tensors are projections of a non-smooth function on the velocity grid. Therefore, they cannot be approximated with high accuracy by a tensor with small ranks.

Nevertheless, as mentioned above, in the formula for the face flux \eqref{eq:exact-flux}, the tensor $|\hat{\xi}_{n,j}|$ can be replaced with some estimate. This can be interpreted as replacing the exact numerical flux with a Rusanov-type flux. In our numerical experiments, we used approximations for $|\hat{\xi}_{n,j}|$ with the Tucker ranks $6$ for all faces. 

\begin{figure}[h]
    \centering
    \includegraphics[width = \textwidth]{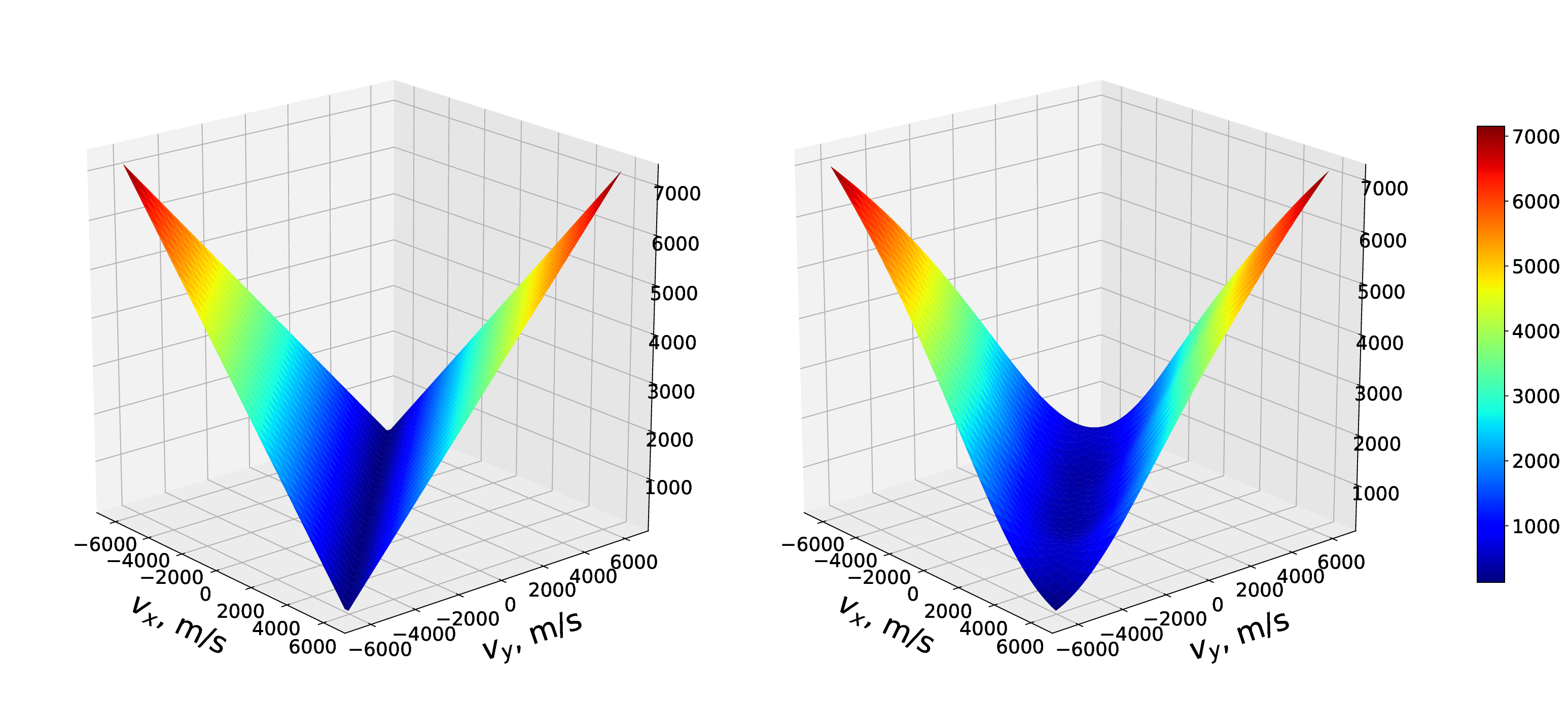}
    \caption{Left: slice of the exact tensor $|\hat{\xi}_{n, j}|(:, :, i_3^0)$, right: slice of the approximation with the Tucker ranks equal to 6.}
    \label{fig: vel-estimate}
\end{figure}

Figure \ref{fig: vel-estimate} shows a comparison between the cross-sections at $i_3 = const$ of the exact tensor $|\hat{\xi}_{n,j}|$ and its low-rank approximation for a random normal vector. It can be seen that the estimate mimics well the exact function. The numerical experiments have shown that the approximation with ranks $6$ is sufficient for the first order scheme. However, for higher-order schemes, a more careful approximation may be needed.

After all the operations in the algorithm which may lead to a large increase in the Tucker ranks, we add rounding with a prescribed relative error $\epsilon$. 
It should be noted that when applying a specific tensor format, it is necessary to consider the computational complexity of each element-wise operation and rounding, and not only the asymptotic growth rate but also the constants involved in the estimates. In the method under consideration, it makes no sense to insert rounding after each operation. In addition, some operations should be reordered. For example, it is preferable to avoid the Hadamard multiplication of two tensors with large ranks, whereas, in contrast, element-wise summations for the same ranks are relatively cheap.

In the implicit scheme, the right-hand side in \eqref{eq:expl-scheme-final} is taken at the $t^{n+1}$ and then linearized:
\begin{align}
\begin{split}
    & \frac{\bm{\hat{f}}^{n+1} - \bm{\hat{f}}^n}{\Delta t} = \bm{\hat{R}}^{n+1} \approx \bm{\hat{R}}^n + \frac{\partial \bm{\hat{R}}^n }{\partial \bm{\hat{f}}^n} (\bm{\hat{f}}^{n+1} - \bm{\hat{f}}^n) \Rightarrow \\
    & \left[ \frac{1}{\Delta t} \hat{I} - \frac{\partial \bm{\hat{R}}^n }{\partial \bm{\hat{f}}^n}\right] (\bm{\hat{f}}^{n+1} - \bm{\hat{f}}^n) = \bm{\hat{R}}^n .
\end{split}
    \label{eq:linear-system}
\end{align}
Here, the bold symbols with hat denote vectors consisting of tensors, for example: 
\[\bm{\hat{f}}^{n} = [\hat{f}_1^{n}, \ldots, \hat{f}_{N_C}^{n}]^T\]
The elements of the jacobian matrix $\frac{\partial \bm{\hat{R}}^n }{\partial \bm{\hat{f}}^n}$ are tensors consisting of element-wise derivatives for values at each point in the velocity grid. $\hat{I}$ is the diagonal matrix with diagonal elements (tensors) equal to $\hat{\mathbbm{1}}$.
In the LU-SGS and Jacobi iterative methods, we need to compute an element-wise division by the diagonal elements $\hat{D}_i$ of the left-hand side matrix in the system \eqref{eq:linear-system}. If we do not take into account the nonlinear dependence of $\hat{f}_i^S$ on $\hat{f}_i$, then the diagonal tensors have the following form:
\begin{align}
    \hat{D}_i = \hat{\mathbbm{1}} \left(\frac{1}{\Delta t} + \nu(\hat{f}_i^n) \right)  + \frac{1}{2 |V_i|} \sum_{l=1}^{N_F(i)} |A_{j(i,l)}| \, \lvert \hat{\xi}_{n,j(i,l)}\rvert .
    \label{eq:diagonal_term}
\end{align}
From this expression it is clear that the elements of the $\hat{D}_i$ are positive.
There is no algorithm for exact component-wise division in the Tucker format. Division can be computed by the  cross-approximation techniques proposed in  \cite{oseledets2008tucker}. In our method, we adopt a simpler approach. Since the linearization in \eqref{eq:linear-system} is already inexact, one can use a simplified formula for the flux to compute $\frac{\partial \bm{\hat{R}}^n }{\partial \bm{\hat{f}}^n}$. We employ the following coarse bound for the $\lvert \hat{\xi}_{n,j}\rvert$:
\begin{align}
    \lvert \hat{\xi}_{n,j}\rvert(i_1, i_2, i_3) \le \sqrt{\xi_{1, i_1}^2 + \xi_{2, i_2}^2 + \xi_{3, i_3}^2}  .
    \label{eq:xin-one-rank-approx}
\end{align}
This bound is then approximated by a tensor with ranks $1$ which is used in \eqref{eq:diagonal_term} to compute the  approximation $\hat{D}_i^{est}$ with ranks $1$. After that, the division by the $\hat{D}_i^{est}$ is computed exactly in $O(n r^3)$ operations:
\begin{small}
\begin{align}
\begin{split}
    & \hat{B}(i_1, i_2, i_3) = \frac{\hat{A}(i_1, i_2, i_3)} {\hat{D}_i^{est}(i_1, i_2, i_3)} = \\[10pt]
    & \frac{\displaystyle \sum_{k_1, k_2, k_3} G^A(k_1, k_2, k_3) U^A_1(i_1, k_1) U^A_2(i_2, k_2) U^A_3(i_3, k_3) }{U_1^D(i_1) U_2^D(i_2) U_3^D(i_3)} = \\[10pt]
    & \sum_{k_1, k_2, k_3} G^A(k_1, k_2, k_3) \frac{U^A_1(i_1, k_1)}{U_1^D(i_1)} \frac{U^A_2(i_2, k_2)}{U_2^D(i_2)} \frac{U^A_3(i_3, k_3)}{U_3^D(i_3)} .
    \end{split}
    \label{eq:one_rank_division}
\end{align}
\end{small}

\begin{figure}
    \centering
    \includegraphics[width = \textwidth]{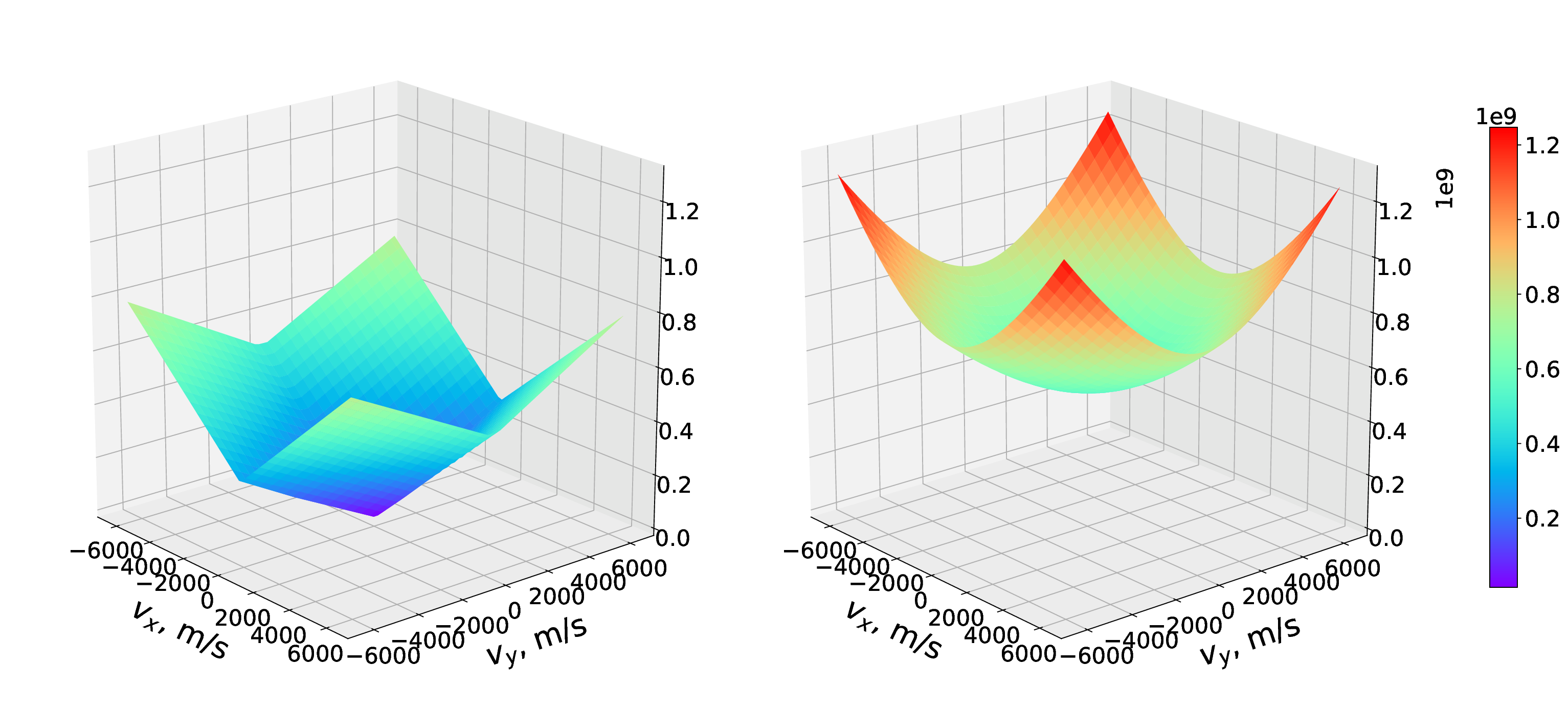}
    \caption{Exact diagonal tensor $\hat{D}_i$ (left) and its approximation with ranks $1$ (right).}
    \label{fig:diags}
\end{figure}

Figure \ref{fig:diags} shows the cross-sections of the $\hat{D}_i$ and the $\hat{D}_i^{est}$ for a space cell. It is clear that the approximation is very rough, but the numerical experiments show that it still provides faster convergence compared to the explicit method, which has a severe stability restriction on the time step $\Delta t$.
\section{Implementation}
\label{sec:impl}
For comparison between the two methods, both standard discrete velocity method and its tensorized version are implemented in the Python language. The program consists of three main Python modules: 
\begin{enumerate}
    \item \textit{read\_starcd.py} -- an auxiliary module for reading an unstructured mesh in the StarCD format. It contains class \textit{Mesh}. The constructor of this class takes path to the folder with mesh files and creates an object where all the information needed in the numerical method is stored (cell volumes, face normals, etc.) This object is  then serialized using the \textit{pickle} module. Afterward, in the run script the mesh object is read from the serialization file. 
    \item \textit{solver.py} -- this module declares a class for the solution and the class method \textit{make\_time\_steps} to perform a number of time steps using the described first order discrete velocity method operating with full tensors.
    \item \textit{solver\_tucker.py} -- contains the implementation of the tensorized version of the discrete velocity method.
\end{enumerate}
Besides, there are two scripts for two test problems: the first is the 1D shock wave structure problem, and the second is the flow around a planar circular cylinder (see section \ref{sec:test}). The shock wave test can be used for the first validation and experiments since the spatial mesh is very small, and so is the computational time. The second test demonstrates that the tensorized version of the algorithm  provides a significant memory reduction in real-life problems.

The spatial mesh for any problem can be created using appropriate software. StarCD is a widespread format, so one can convert a mesh from other formats to the StarCD format. We used Ansys ICEM to create meshes for our tests.

To solve a new problem, one needs to create an object of the \textit{Problem} class (see listing \ref{list:problem}) and pass it to the constructor of the object of the \textit{Solution} class (listing \ref{list:solution}) together with objects of the \textit{Mesh} class and the \textit{VelocityGrid} class, which contains tensors for the velocity grid in the Tucker format.
\begin{lstlisting}[language=Python, caption=\textit{Problem} class, basicstyle = \small, label = list:problem, frame=single]
class Problem:
    def __init__(self, bc_type_list = None, 
                 bc_data = None, f_init = None):
        # list of boundary conditions' types
        # according to order in starcd '.bnd' file
        # list of strings
        self.bc_type_list = bc_type_list
        # data for b.c.
        # list of lists
        self.bc_data = bc_data
        # Function to set initial condition
        self.f_init = f_init
\end{lstlisting}
For example, in listing \ref{list:bc}, the  boundary and the initial condition for the flow past cylinder is defined. Currently, a basic set of boundary conditions is implemented, including in-out conditions (which are actually the same free stream condition), wall boundary condition \eqref{eq:wall-bc-f} and symmetry in each coordinate direction.
\begin{lstlisting}[language=Python, caption=Setting initial and boundary condition for the flow past cylinder, basicstyle = \small, label = list:bc, frame = single]
f_init = lambda x, y, z, v: tuck.tensor(
        solver.f_maxwell_t(
            v, n_in, u_in, 0., 0., T_in, gp.Rg))
    
f_bound = tuck.tensor(
        solver.f_maxwell_t(
            v, n_in, u_in, 0., 0., T_in, gp.Rg))

fmax = tuck.tensor(
        solver.f_maxwell_t(
                v, 1., 0., 0., 0., T_w, gp.Rg))
                
problem = solver.Problem(
        bc_type_list = 
        ['sym-z', 'in', 'out', 'wall', 'sym-y'],
        bc_data = [[],
                   [f_bound],
                   [f_bound],
                   [fmax],
                   []], f_init = f_init)
\end{lstlisting}
\begin{minipage}{\linewidth}
\begin{lstlisting}[language=Python, caption=\textit{Solution} class, basicstyle = \small, label = list:solution, frame = single]
class Solution:
    def __init__(self, gas_params, problem, 
                                    mesh, v, config):
        # initialize all required tensors and initial
        # values of the solution
        ...
        if (config.init_type == 'default'):
            for i in range(mesh.nc):
                x = mesh.cell_center_coo[i, 0]
                y = mesh.cell_center_coo[i, 1]
                z = mesh.cell_center_coo[i, 2]
                self.f[i] = problem.f_init(x, y, z, v)
        elif (config.init_type == 'restart'):
            # restart from distribution function
            self.f = self.load_restart()
            
    def make_time_steps(self, config, nt):
        # perform nt time steps
        ...
\end{lstlisting}
\end{minipage}
\section{Test problem}
\label{sec:test}
\begin{figure}[h!]
    \centering
    \includegraphics[width=0.8\textwidth]{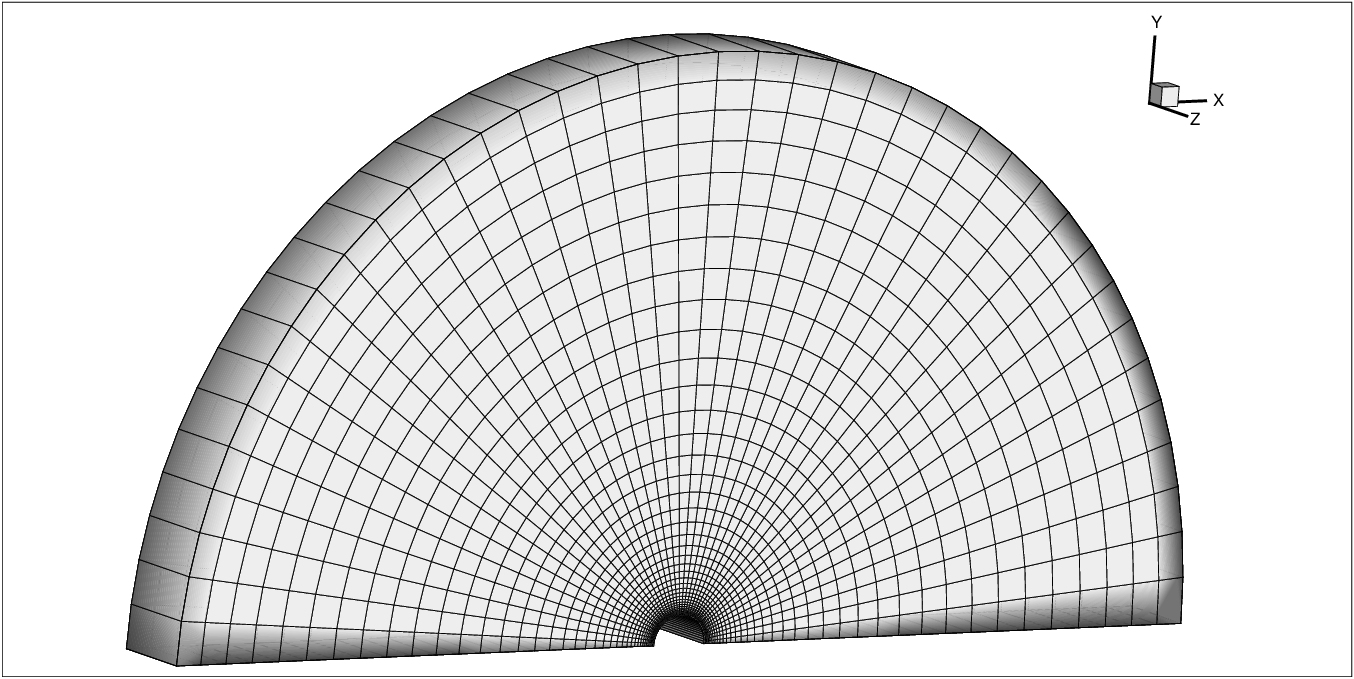}
    \caption{Computational domain and mesh for the test problem.}
    \label{fig:mesh-cyl}
\end{figure}

The problem of high-speed rarefied gas flow past a circular cylinder is considered.
The setup of the problem is taken from \cite{Lofthouse:2008a}. The solutions by the S-model equation  and the exact Boltzmann equation was compared against the DSMC solution in several recent papers
\cite{titarev2018application,Frolova:2018b,Titarev:2018c}
for large free-stream Mach numbers (up to 25) and good agreement was observed. 
The geometry of the computational domain, along with the spatial mesh, is shown in  figure \ref{fig:mesh-cyl}. The problem is essentially two-dimensional, but we solve it as a 3D problem on the 3D spatial mesh with one cell along the $z$-axis. The hexahedral mesh in the physical space is treated as unstructured by the solver. The free stream flow is directed along the $x$-axis. The boundary condition \eqref{eq:wall-bc-f} is set on the wall. 
At the outer boundary, the free stream condition is set. At the remaining boundaries, the symmetry boundary condition is used.

The following dimensional parameters were chosen for the free stream:
$v_0 = 2630$ m/s, $ n_0 = 2 \cdot 10^{23} \mbox{ m}^{-3}$,  $T_0 = 200 \, K$. The wall temperature $T_w = 5 T_0$, the  cylinder radius $r = 1.35 \cdot 10^{-5}\,$ m. The Knudsen number calculated by the parameters of the free stream and the radius of the cylinder $Kn \approx 0.56$, the  Mach number equals to $10$. The power law  was used for viscosity:
\begin{align}
    \mu(T) = \mu_0 \left(\frac{T}{T_0}\right)^{0.734}, \; \mu_0 = 1.61 \cdot 10^{-5} \; \mbox{Pa} \cdot \mbox{s}, T_0 = 200 \, \mbox{K} .
\end{align}

The uniform velocity grid in the cube $[-\xi_{max}, \xi_{max}]^3$ with $\xi_{max} \approx 6400 \,\mbox{m/s}$ contains $N = 64$ nodes in each direction. The number of cells in the spatial mesh equals to $1600$. For this test case, we choose a relatively coarse spatial mesh so that the baseline method can be run on a desktop computer with a rather low RAM (storage of all the values of the distribution function requires more than 3 GB of memory). Therefore, near the surface, the mesh resolution is poor (the normal size of the first cell is too large), but  we concentrate on comparing the two methods rather than accurate computation of the wall friction and heat transfer.

\begin{figure}[h!]
    \centering
    \includegraphics[width = 0.8\textwidth]{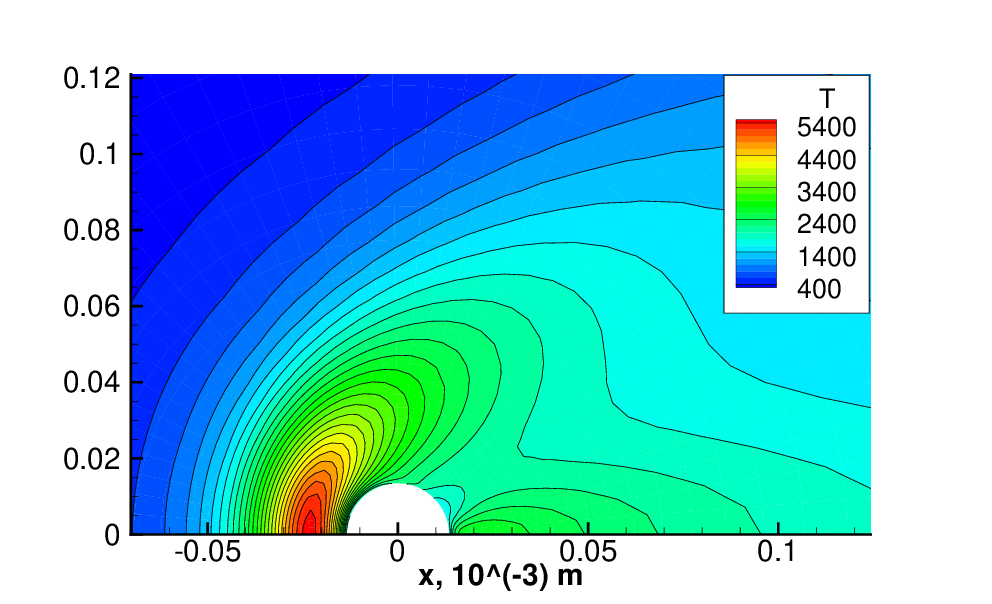}\\
    \includegraphics[width = 0.8\textwidth]{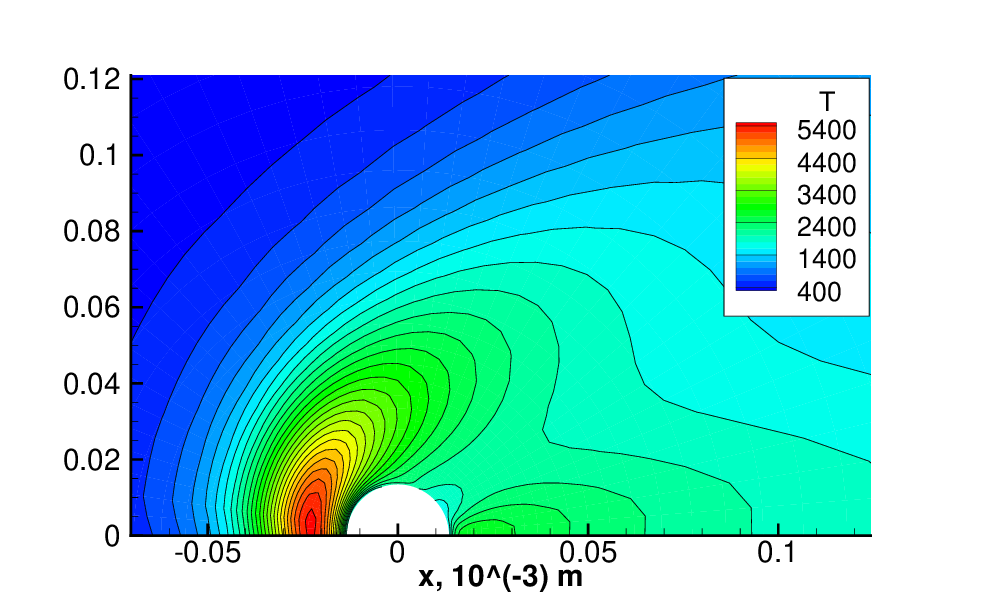}
    \caption{Temperature distribution. Top: standard method, bottom: tensorized method ($\epsilon = 10^{-5}$).}
    \label{fig:temp-color}
\end{figure}

\begin{figure}[h!]
    \centering
    \includegraphics[width = 0.45\textwidth]{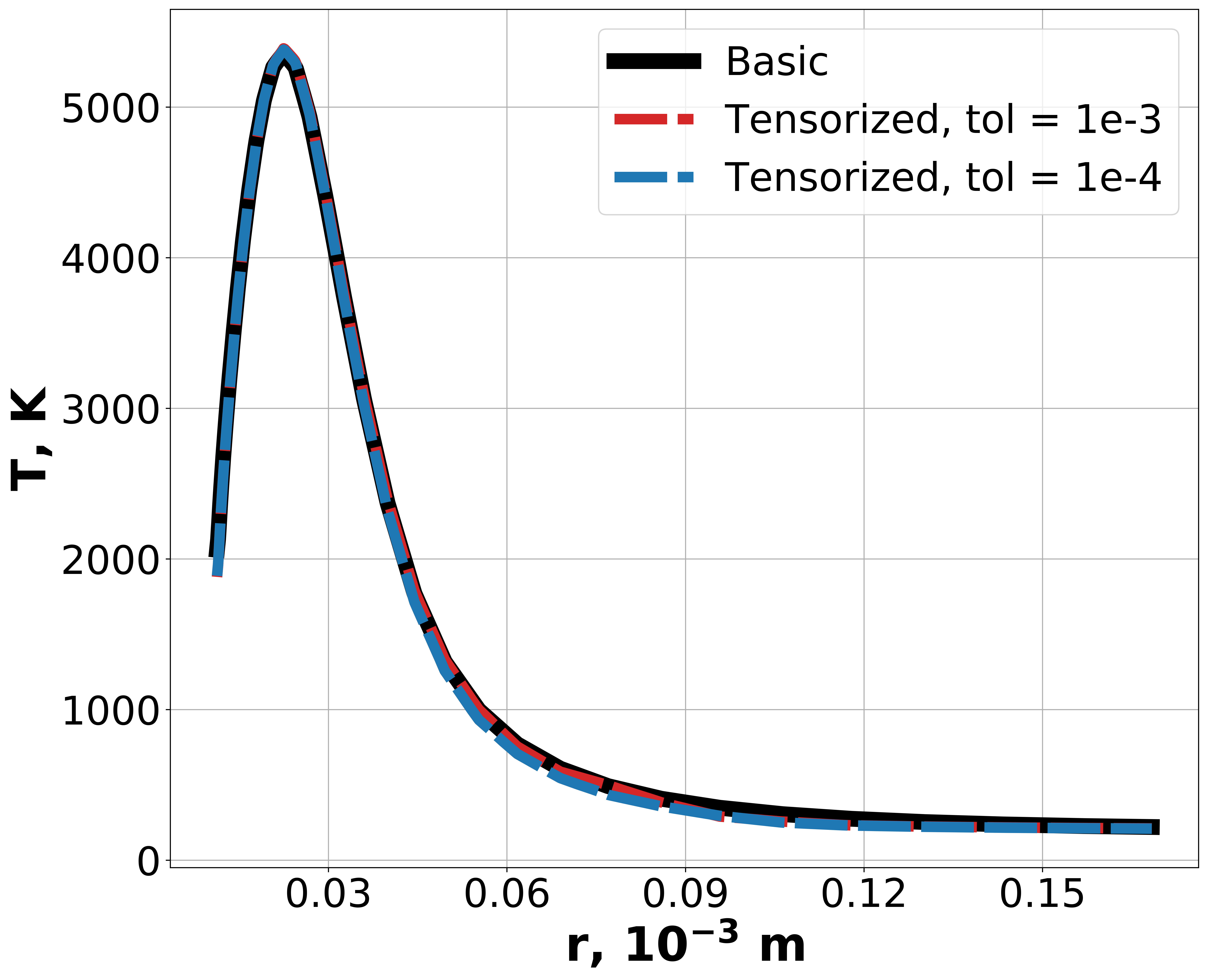} 
    \includegraphics[width = 0.45\textwidth]{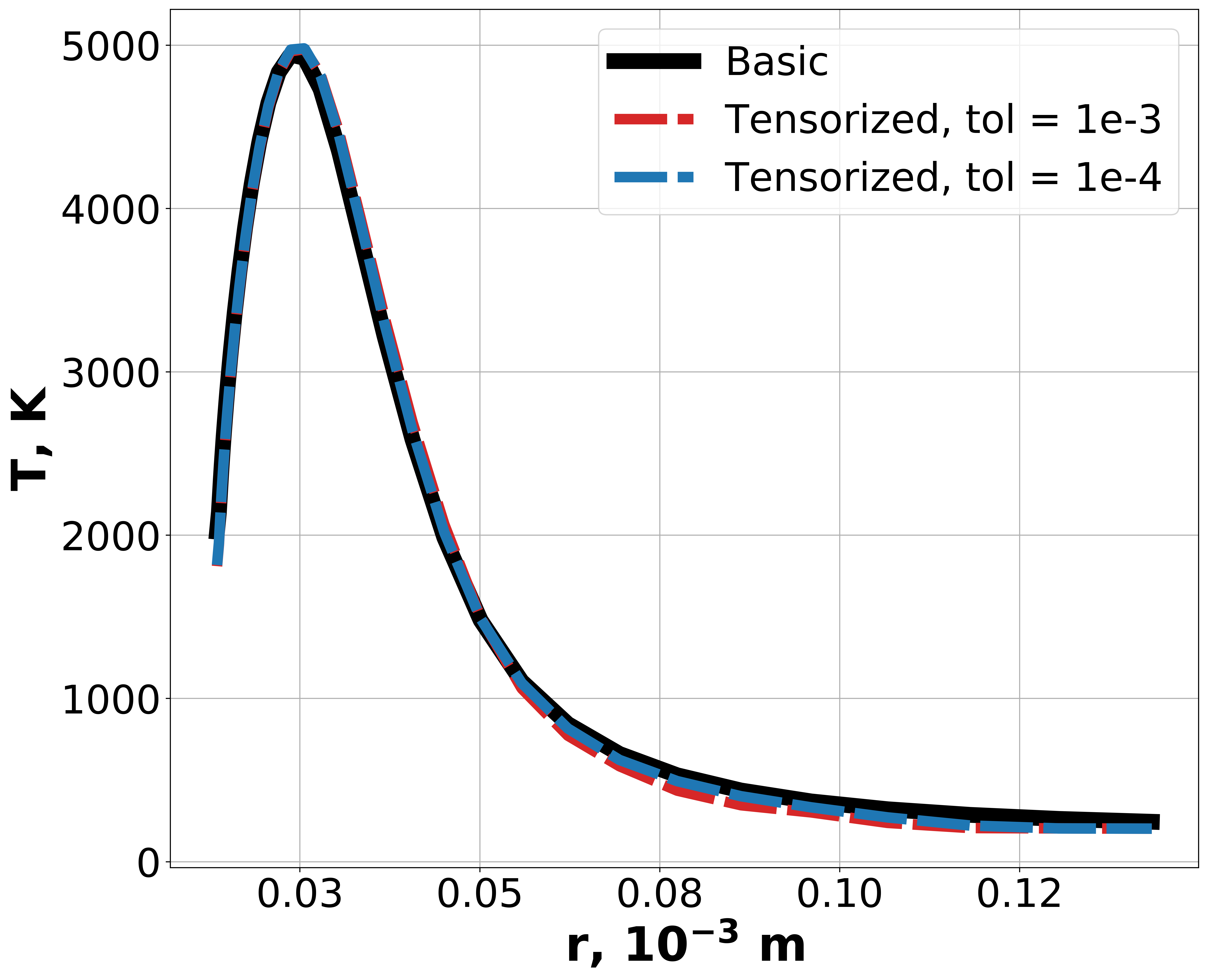} 
    \caption{Temperature profiles along the stagnation line (left) and along the normal at 35 degrees (right).}
    \label{fig:temp-1d}
\end{figure}

Figure \ref{fig:temp-color} shows the temperature distribution obtained  by the standard and tensorized methods. A typical structure of the flow with a detached shock wave can bee seen. Figure \ref{fig:temp-1d} shows the temperature profiles versus the normal coordinate for the stagnation line and for the line normal to the cylinder at an angle of 35 degrees from the stagnation line. The temperature was chosen  for comparison since it is a more sensitive quantity, and the differences for density and velocity are much smaller. 

\begin{figure}[h!]
    \centering
    \includegraphics[width = \textwidth]{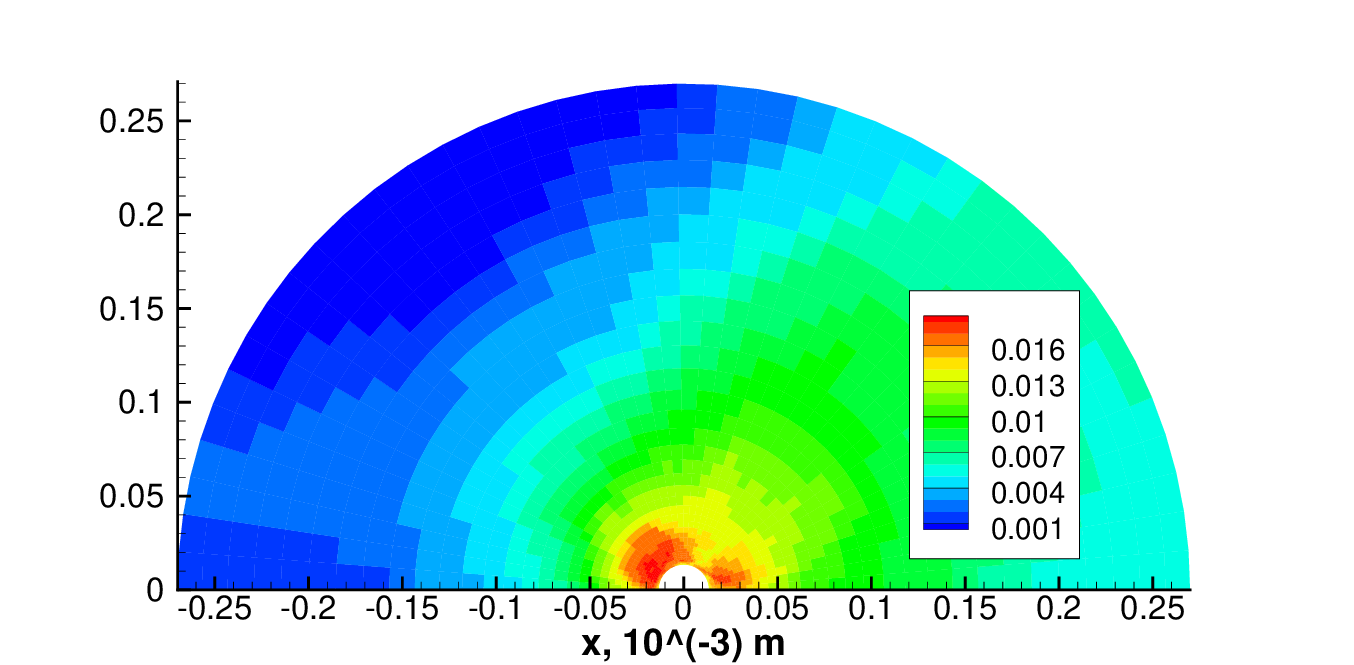}
    \caption{Distribution of the compression ratio $C$ for $\epsilon = 10^{-5}$.}
    \label{fig:cyl-ranks}
\end{figure}

Figure \ref{fig:cyl-ranks} demonstrates the distribution of the compression ratio 
\begin{align*}
C = (r_1 r_2 r_3 + N(r_1 + r_2 + r_3)) / N^3 ,
\end{align*}
where $r_\alpha$ are the Tucker ranks of the tensor $\hat{f}$ in each cell. It can be seen that the rounding of the Tucker tensors works like an adaptive mesh refinement: near the inflow, where the distribution function is almost equilibrium, the ranks are very small. Near the shock wave and the surface, the ranks are automatically increased  to provide the prescribed accuracy.

We computed several solutions using the tensorized method with different values of the relative accuracy in the rounding. Table \ref{tab:table} shows how the relative accuracy $\epsilon$ affects the total compression ratio (averaged over all the cells) and the relative difference in the $2$-norm between temperature distributions in the baseline and tensorized solutions. The storage for the solution in the tensorized method is about $100$ times smaller than in the standard method even for the smallest $\epsilon = 10^{-5}$. The results show that the accuracy $\epsilon = 10^{-4}$ is sufficient for the tensorized method, because the solution does not change significantly with further decrease in $\epsilon$. Moreover, it is clear that the solution by the tensorized method does not converge to the baseline solution. The reason is that a different numerical method is actually used, because we replace the $|\hat{\xi}_{n,j}|$ in \eqref{eq:exact-flux} with an approximation with the Tucker ranks independent of the $\epsilon$. 
\begin{table}[h]
    \centering
    \begin{tabular}{|c|c|c|}
    \hline
    $\epsilon$ & Compression & $\Delta T$  \\
    \hline
    $10^{-3}$ & $5.4 \cdot 10^{-3}$ & $3.3 \cdot 10^{-2}$ \\
    \hline
    $10^{-4}$ & $7.2 \cdot 10^{-3}$ & $2.7 \cdot 10^{-2}$ \\
    \hline
    $10^{-5}$ & $1.1 \cdot 10^{-2}$ & $3.1 \cdot 10^{-2}$ \\
    \hline
    \end{tabular}
    \caption{Dependence of the compression ratio and the relative temperature difference on the relative accuracy of the rounding.}
    \label{tab:table}
\end{table}

\begin{figure}[h!]
    \centering
    \includegraphics[width = 0.9\textwidth]{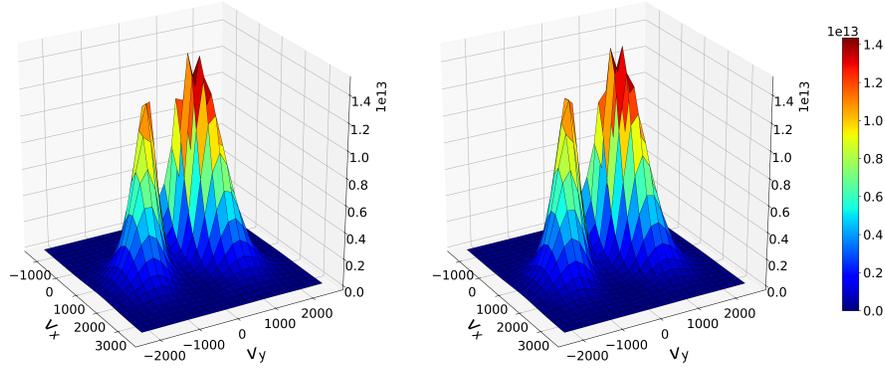}
    \caption{Slice of the distribution function tensor in one space cell. Left - standard method, right - tensorized method with $\epsilon = 10^{-5}$.}
    \label{fig:distr-function}
\end{figure}

Figure \ref{fig:distr-function} shows the $z$-slice of the distribution function in the cell with $x = 2.46 \times 10^{-5}, y = 10^{-6}$ (near the stagnation line on the shock front). In this area, the flow is strongly non-equilibrium, and the distribution function has two peaks. It can be seen that the difference between two solutions is negligible, i.e., the  tensorized method successfully captures the main properties of the distribution function.

Despite the significant memory reduction, the computational time of both methods is approximately equal for this test. For low values of the $\epsilon$, the tensorized method is slower then the baseline method. The reason is the high cost of the element-wise product and rounding. The same situation is reported in other studies, for instance \cite{Dolgov2014268}. However, for this test we used a rather small velocity grid ($64^3$ nodes). For larger grids, the tensorized algorithm would be faster than the standard one.

\section{Concluding remarks and perspectives}
\label{sec:conclusion}
The Boltzmann-T solver for the numerical solution of the kinetic Boltzmann equation with a model collision integral is described. The solver provides a working example of the implementation of a tensorized discrete velocity method. This implementation demonstrates prospects of using tensor decompositions for significant memory reduction in practical computations by the discrete velocity method on an unstructured spatial mesh.

We draw the following conclusions from our experience, which may be useful for other researchers dealing with tensorized versions of a discrete velocity method: 
\begin{enumerate}
    \item Problem with tensors generated by a non-smooth function (such as~$|\hat{\xi}_{n,j}|$) can be overcome if the spherical coordinate system is used in the velocity space. In this case,  tensors like~$|\hat{\xi}_{n,j}|$ have low ranks. One possible drawback is that spherical coordinates may lead to more complicated quadrature formulas.
    \item In this paper, we consider the most straightforward approach for the algorithm modification: all basic operations are replaced with tensorized analogs. The more elegant approach is to use cross-approximation techniques and it should be a subject of a future research. Nevertheless, in our opinion, the straightforward approach is very robust and does not require a deep understanding of underlying tensor algorithms.
    \item In all standard tensor formats, storage and complexities of the main algorithms are proportional to $n$ -- the size of the original tensor in one dimension. For the large $n$, artificial increase of the dimensionality or the so-called quantized tensor formats \cite{dolgov2013QTT} can be applied to decrease memory consumption even further.
\end{enumerate}

In future, we plan to implement a parallel version of our solver using \textit{mpi4py} package and spatial mesh decomposition. Besides, we plan to add model collision integrals for diatomic gas with internal degrees of freedom. The numerical method will be extended to higher orders of approximation, tetrahedral spatial meshes, and unsteady problems.

\section*{Acknowledgements}
The authors thank Sergey Dolgov, Maxim Rakhuba, and Ivan Oseledets for their helpful recommendations regarding tensor formats.

Some of the calculations were run on supercomputers of the Joint Supercomputing Center of the Russian Academy of Sciences.

A. Chikitkin and E. Kornev were supported by the Presidential Grant  number MK-2855.2019.1.





\bibliographystyle{elsarticle-num}
\bibliography{references}







\end{document}